\documentstyle[12pt]{article}

\makeatletter\@addtoreset{equation}{section}\makeatother

\newtheorem{Tm}{Theorem}[section]
\newtheorem{Rk}{Remark}[section]
\newtheorem{Lm}{Lemma}[section]
\newtheorem{Co}{Corollary}[section]
\newtheorem{Df}{Definition}[section]
\newtheorem{Pn}{Proposition}[section]

\begin{document}
\begin{center}
\begin{Large}
{\bf Laguerre Entire Functions and Related Locally Convex Spaces  }
\end{Large}
\vskip.8cm
{\it by}
\vskip.8cm
\begin{large}
Yuri Kozitsky \footnote{Supported in part under the
Grant KBN No 2 P03A 02915} \ \ \ \ Lech Wo{\l}owski
\end{large}
\end{center}

\vskip.8cm

\def\cit#1{\cite{[#1]}}
\def\kasten{\hfil\vrule height6pt width5pt depth-1pt\par }
\def\C{{C\!\!\! C}}
\def\R{{I\!\! R}}
\def\N{{I\!\! N}}
\def\Z{{Z\!\!\! Z}}
\def\pp{\partial}
\def\ds{\displaystyle}
\def\half{{1 \over 2}}
\def\lra{\longrightarrow}
\def\Ra{\Rightarrow}
\def\a{\alpha}
\def\b{\beta}
\def\g{\gamma}
\def\d{\delta}
\def\e{\epsilon}
\def\i{\iota}
\def\k{\kappa}
\def\t{\tau}
\def\l{\lambda}
\def\om{\omega}
\def\s{\sigma}
\def\vp{\varphi}
\def\ve{\varepsilon}
\def\vt{\vartheta}
\def\D{\Delta}
\def\G{\Gamma}
\def\L{\Lambda}
\def\Om{\Omega}
\def\AC{{\cal A}}
\def\Aa{{\cal A}_a }
\def\BC{{\cal B}}
\def\CC{{\cal C}}
\def\DC{{\cal D}}
\def\EC{{\cal E}}
\def\FC{{\cal F}}
\def\GC{{\cal G}}
\def\HC{{\cal H}}
\def\IC{{\cal I}}
\def\KC{{\cal K}}
\def\LC{{\cal L}}
\def\Lp{{\cal L}^{+}}
\def\LM{{\cal L}^{-}}
\def\Ln{{\cal L}_0 }
\def\La{{\cal L}_a }
\def\Lpa{{\cal L}^{+}_{a}}
\def\NC{{\cal N}}
\def\PC{{\cal P}}
\def\XC{{\cal X}}
\def\lv{\left\vert}
\def\rv{\right\vert}
\def\bb#1{ {\lv #1 \rv} }
\def\bq#1{ {\lv #1 \rv}^2 }
\def\BE#1{ {\left\Vert #1 \right\Vert} }
\def\ti{\; \times \!\!\!\!\!\!\!\! \mathop{\phantom\sum}}
\def\dop{\dot +}
\def\C{\hbox{\vrule width 0.6pt height 6pt depth 0pt \hskip -3.5pt}C}
\def\ope{\;
\lra_{\!\!\!\!\!\!\!\!\!\!\!\!\!_{\hbox{$_{n\to\infty}$}}} \;}
\def\be{\begin{equation}}
\def\th{\theta}
\def\ee{\end{equation}}
\def\SC{\cal S}
\def\beq{\begin{eqnarray}}
\def\eeq{\end{eqnarray}}
\def\k{\kappa}
\def\ra{\rightarrow}
\def\bk{\bar{\kappa}}
\def\z{\zeta}
\def\CN{{\cal N}}
\def\sta{{\stackrel{\rm def}{=}}}
\def\Dth{{\Delta}_{theta}}
\def\TC{{\cal T}_C}
\pagenumbering{arabic}
\begin{abstract}
A scale $\{{\AC}_a, \ a\geq 0\}$, ${\AC}_a \subset
{\AC}_b $, for $a\leq b$,
of the Fr{\'e}chet spaces of exponential type entire functions
 of one complex variable is considered. Certain
special properties of the subsets
of ${\AC}_a$ consisting of Laguerre entire functions, which are
obtained as uniform limits on compact subsets of $\ \C$ of
polynomials with real nonpositive zeros only, are described.
On the space ${\AC}_b$, the operators having the form
$\vp (\Delta_\theta )$, where
$\vp \in {\AC}_a$, $ab<1$, and $\Delta_\theta = (\theta + zD)D$
with $\theta \geq 0$ and $D=d/dz $, are defined. They are shown to
preserve the set of Laguerre entire functions $\Lp$. An integral form
of $\exp (a \Delta_\theta )$ with $a>0$ is found
that allows to construct some extensions of this operator.
These results are used to obtain and to study the solutions of a
certain initial value problem involving $\Delta_\theta$.
\end{abstract}

\section{Introduction and Main Results}

\subsection{Introduction}

In this paper, topological vector spaces of exponential type entire
functions of one complex variable are considered. We introduce a
scale of spaces $\{{\AC}_a, \ a\geq 0\}$, ${\AC}_a \subset
{\AC}_b $, for $a\leq b$, where each ${\AC}_a$ is defined as
a Fr{\'e}chet space. Applying some results of \cit{Tay} we show
that the set of all polynomials is dense in every ${\AC}_a$ and that
the relative topology on every bounded subset $B\subset{\AC}_a$
coincides with the topology of uniform convergence on compact
subsets of $\ \C$ (Theorem \ref {n1tm}). This implies that
each $\Aa $ inherits Montel's property from the space of all
entire functions, i.e., its bounded subsets
are precisely the relatively compact subsets.
It also means if
a sequence $\{f_n \}\subset {\AC}_a$ (a) converges
uniformly on compact subsets of $\ \C$ to a function $f$;
(b) this $f$ belongs to $ {\AC}_a$; (c) it is bounded in ${\AC}_a$,
then this sequence converges to $f$ also in ${\AC}_a$.
 We study the possibility to slack
the above conditions by excluding (c) and to obtain it as a result
of the convergence $f_n \ra f$ in ${\AC}_a$. It turns out that this
generally is impossible and we give some examples of sequences
$\{f_n \}\in {\AC}_a$ converging to $f\in{\AC}_a$ uniformly on
compact subsets of $\ \C$, which are unbounded in ${\AC}_a$.
At the same time we prove that this possibility occurs
(Theorem \ref{n2tm})
for the sequences $\{f_n \}  \subset {\Lp}$. The latter is
a class of Laguerre entire functions.
These functions are
obtained as uniform limits on compact subsets of $ \ \C$ of the
sequences of polynomials possessing real nonpositive zeros only.
The class $\Lp$ was being studied by many authors
during all this century, that was caused by a number of significant
properties and various applications of these functions. A considerable
survey of this study is given in \cit{Il} (see also \cit{Lev}).

On the introduced Fr{\'e}chet spaces ${\AC}_b$, we define
the operators having the form $\vp (\Delta_\theta )$, where
$\vp \in {\AC}_a$, $ab<1$, and $\Delta_\theta = (\theta + zD)D$
with $\theta \geq 0$ and $D=d/dz $. It is proven that
each such an operator continuously maps
${\AC}_b \ra {\AC}_c $ with $c =b(1-ab)^{-1}$. The main result
obtained here is Theorem \ref{a1tm} which asserts that if $\vp
\in {\Lp}$ and $f\in {\Lp}$, then $\vp(\Delta_\theta )f\in {\Lp}$
provided it exists in some $ {\AC}_c$. To prove this assertion
we construct a technique in the form of Lemma \ref{a1lm}  --
Lemma \ref{a4lm},
which allows us to control the distribution of zeros of the
function $\vp(\Delta_\theta )f$. Further we find an integral form
of $\exp(a\Delta_\theta )$, $a\geq 0$ (Proposition
 \ref{a2pn}), which together with an analog of the operation rules
(Proposition \ref {a3pn}) allows to construct some
extensions of this operator. This is then used to obtain and to
study the solutions of a certain initial value problem involving
$\Delta_\theta $ (Theorem \ref{e1tm}).

All simple proofs follow directly
the statements. More complicated and more technical proofs are placed
in the second section.

\subsection{Spaces of Entire Functions}

Let $\EC$ be the set of all entire functions $ \ \C \ra \C$
equipped with the pointwise linear operations and with
the topology $\TC $ of uniform convergence on compact subsets of $\ \C$.
 For $b>0$, we define
$$
{\cal B}_{b}=\{f\in {\EC }\mid \Vert f\Vert _{b}<\infty \};
$$
where
\begin{equation}
\label{1}
\Vert f\Vert _{b}=\sup_{k\in \N_0 } \{ b^{-k}\mid f^{(k)}(0)\mid   \};
\ \ \  f^{(k)}(0)=(D^{k}f)(0) = \frac{d^k f}{dz^k } (0),
\end{equation}
and $\N_0 $ stands for the set of all nonnegative integers.

\begin{Pn}
\label{1pn}
Every $\left( {\cal B}_{b},\Vert \cdot \Vert _{b}\right) $ is a Banach
space. Every sequence $\{f_n \}$ converging in ${\BC}_b $, converges also in
$({\EC} , \TC ) $.
\end{Pn}

{\bf Proof.} To prove the first part of this statement we need
only to show
the completeness of ${\BC}_b $.
Let $\{f_{n}\}$ be a Cauchy sequence in ${\cal B}_{b}$,
$K$ be a compact subset of $ \ \C$, and
$$
r(K) \ \sta \ \sup_{z \in K} \vert z \vert .
$$
Given  $\varepsilon
>0$, we choose $N$ such that for all $m,p>N,$%
\begin{equation}
\label{2}
\left\| f_{m}-f_{p}\right\| _{b}<\varepsilon ,
\end{equation}
and find
\beq
\sup_{z\in K}\left| f_{m}(z)-f_{p}(z)\right| & = & \sup_{z\in K}\left|
\sum_{k=0}^{\infty }\frac{f_{m}^{(k)}(0)-f_{p}^{(k)}(0)}{k!b^{k}}%
(bz)^{k}\right| \nonumber \\
& \leq & \left\| f_{m}-f_{p}\right\| _{b}\sup_{z\in K}\sum_{k=0}^{\infty }\frac{%
\left| bz\right| ^{k}}{k!}\leq \varepsilon \exp(br(K)) , \nonumber
\eeq
 which means that $\{f_{n}\}$ is
a Cauchy sequence in $({\EC},\TC )$. Hence there exists an entire
function $f$ such that $f_{n}\rightarrow f$ in $({\EC}, \TC )$.
This proves the second part of the statement. The convergence
just established means that,
for every  $k\in \N $, the sequence
$\{ f_{n}^{(k)} (0) \} $ converges to $ f^{(k)} (0)$
(by Weierstrass' theorem).
On the other hand, for every $ k\in \N_0 $,
one easily gets from (\ref{1}), (\ref{2})
\begin{equation}
\label{3}
b^{-k}\left| f_{m}^{(k)}(0)-f_{p}^{(k)}(0)\right| \leq \left\|
f_{m}-f_{p}\right\| _{b}<\varepsilon .
\end{equation}
Thus passing here to the limit $p\rightarrow \infty $ one obtains
\be
\label{4}
\sup_{k\in \N_{0}}\left\{ b^{-k}\left|
f_{m}^{(k)}(0) -f^{(k)}(0) \right| \right\} \leq \varepsilon .
\ee
Hence
\beq
\sup_{k\in \N_{0}}\left\{ b^{-k}\left|f^{(k)}(0) \right| \right\} & \leq &
\sup_{k\in \N_{0}}\left\{b^{-k}\left|f^{(k)}_m (0) \right| \right\} +
\sup_{k\in \N_{0}}\left\{b^{-k}\left|f^{(k)}_m (0) - f^{(k)} (0) \right| \right\} \nonumber \\
& \leq &
\|f_m \|_b + \varepsilon , \nonumber
\eeq
which means $f\in \BC_b $. Then the estimate (\ref{4})
implies
$$
\|f_m - f \|_b \leq \ve < 2\ve ,
$$
thus $%
f_{m}\rightarrow f$ in ${\cal B}_{b} $.
\kasten

The central role in this work is played
by the following space of entire functions.
For $a\geq 0$, let
\begin{equation}
\label{5}
{\cal A}_{a}=\bigcap_{b>a}{\cal B}_{b}=\{f\in {\EC}
\mid (\forall b>a ) \ \Vert f\Vert _{b}<\infty \},
\end{equation}
This set equipped with  the topology $%
{\cal T}_{a}$ defined by the family $\{\Vert .\Vert _{b}, \
b>a\}$ becomes a separable locally convex space.
Due to the
following obvious inequality
$$
\Vert f\Vert _{b^{\prime }}<\Vert f\Vert _{b},\qquad b^{\prime }>b
$$
the topology ${\cal T}_{a}$ may also be defined by a countable family of
norms $\{\Vert .\Vert _{b_{n}},n\in \N \}$ (e.g. $%
b_{n}=a+1/n).$ Thus the topological vector space
$({\cal A}_{a},{\cal T}_{a})$, being the
projective limit of a countable family of the Banach spaces
$({\cal B}_{b_{n}},\Vert .\Vert _{b_{n}}), $
is complete and metrizable
(see \cit{Sch}, chap.I, \S\ 6 and
chap.II, \S\ 5). Therefore
$({\cal A}_{a},{\cal T}_{a})$ is a Fr\'{e}chet
space. In the sequel we will write $\EC$, ${\cal A}_{a}$, and
${\cal B}_b$ instead of $ ({\EC}, {\cal T}_C )$, $({\cal A}%
_{a}, {\cal T }_{a})$, and $({\cal B}_b, \BE{.}_b)$ respectively,
assuming that the mentioned topologies
are the standard ones on
these sets.
It should be pointed out that $\EC$ and ${\AC}_0 $, equipped
also with the pointwise multiplication, become algebras.
As a subset of $\EC$, $\Aa$ equipped with ${\cal T}_a $ should
inherit some properties of $(\EC, {\cal T}_C)$. An important
example here is Montel's property -- the bounded subsets of
$\EC$ are precisely the relatively compact subsets (see
\cit{BG} p. 120). A subset $B\subset \EC $ is said to be bounded
if for every compact $K\in \ \C$, there exists $C(K)$ such that
$$
\sup_{f\in B}\{\sup_{z\in K} \lv f(z) \rv \} \leq C(K).
$$
\begin{Df}
\label{Bound}
A subset $B\subset \AC_a $ is said to be bounded in $\Aa$
if for every $b>a$, there exists a constant $C_b $ such that
$$
\sup_{f\in B}\BE{f}_b \leq C_b .
$$
\end{Df}
\begin{Tm}
\label{n1tm}
For every $a\geq 0$, the space ${\AC}_a $ possesses the properties:
\vskip.1cm
\begin{tabular}{ll}
{\rm (i)}  & the set of all polynomials ${\cal P} \subset \EC$
         is dense in $\AC_a $;\\ [.1cm]
{\rm (ii)}  & the relative topology on every bounded subset
            $B\subset \AC_a $
               coincides \\
             & with the topology of uniform convergence on
                  compact \\
              & subsets of $  \ \C $; \\ [.1cm]
{\rm (iii)} & for every $f\in {\AC}_a$ and $g\in {\AC}_b$, their
             product $fg$ belongs to ${\AC}_{a+b} $.
\end{tabular}
\end{Tm}
\begin{Co}
\label{pMonco}
A subset $B\in \Aa $, which is closed and bounded in $\EC$,
is closed in $\Aa$ provided it is bounded there.
\end{Co}
\begin{Co}{\rm [Montel's Property]}
\label{Monco}
The bounded subsets of $\Aa$ are exactly the relatively compact
subsets.
\end{Co}
\begin{Co}
\label{new1co}
Let $\{f_n , \ n\in \N \}$ converge in ${\AC}_a $ to $f$ and
$\{g_n , \ n\in \N \}$ converge in ${\AC}_b $ to $g$. Then the
sequence $\{f_n g_n , \ n\in \N \}$ converges to $fg$ in ${\AC}_{a+b} $.
\end{Co}

It may seem that a sequence $\{f_n , n \in\N\} \subset {\AC}_a$
with the
following properties:
(a) it converges in $\EC$ to a function $f$;
(b) $f \in{\AC}_a $; should possess the
property: (c) it converges to this $f$ also in ${\AC}_a $.
In view
of claim (ii) of the above theorem,
to prove this conjecture it would be enough to show
the boundedness of this sequence in ${\AC}_a $.
In fact, the conjecture is false. Furthermore, a claim like (ii) does not
hold for the Banach spaces ${\cal B}_b$. The following examples
show the subtlety of the situation with the topologies of the
mentioned spaces. Consider the sequence $\{f_n (z) = z^n /n! , \ n
\in \N \}$. Since it consists of polynomials, it is a subset of
${\AC}_0$ and of all ${\cal B}_b $, $b>0$.
By means of (\ref{1}), one finds $\BE{f_n}_a = a^{-n}$.  Thus
the sequence is bounded only in ${\AC}_a$ and in ${\cal B}_a $ with
$a\geq 1$. On the other
hand, this sequence converges in $\EC$ to the function
$f(z) \equiv 0$ which also belongs to all mentioned spaces.
In all ${\AC}_a$, $a\geq 1$, where it is bounded,
the sequence converges to $f\equiv 0$, as it
is prescribed by claim (ii). In ${\AC}_a $, $a<1$ our sequence
possesses (a) and (b) but is unbounded in ${\AC}_a $
hence does not possess (c).
 Moreover,
it is bounded in
${\cal B}_1 $ and does not converge
in this space, i.e., a claim like (ii) would fail
for the Banach spaces ${\cal B}_b $.
Another example, which is only a slight variation of the previous one, is
as follows. Let $f$ be an arbitrary element of ${\AC}_a $, $a<1$.
For
$\ve >0$, we define the sequence of functions $\{f_n (z,\ve )\}$ by
their derivatives at $z=0$:
$$
f_n^{(k)} (0, \ve ) =  f^{(k)} (0) + \ve \delta_{kn},
  \ \ k,n \in \N_0 ,
$$
where $\d$ stands for the Kronecker $\d$--symbol.
Then this sequence converges to $f$ in every ${\AC}_a $, $a\geq 1$,
and hence does in $\EC$. But, for every positive $\ve$,
 it is unbounded in any ${\AC}_a $ with
$a<1$.

 Nevertheless, there exists the subset of $\EC$
such that on its intersections with the spaces ${\AC}_a $, $a\geq 0$,
the above
mentioned properties (a) and (b) imply (c).

\begin{Df}
\label{1df}
A family ${\cal L}$ (respectively ${\cal L}_0 $, ${\cal L}^{+}$,
${\cal L}^{-}$) consists of all entire functions
 possessing the following representation
\begin{equation}
\label{6}
f(z)=Cz^{l}\exp (\a z)\prod_{j=1}^{\infty }(1+\b_{j}z);
\end{equation}
$$
C\in \ \C;\ l\in \N_{0};\ \ \b_{j}\geq \b_{j+1} \geq 0;
\ \sum_{j=1}^{\infty }\b_{j}<\infty ,
$$
with $\a\in \R$  (respectively $\a=0$, $\a\geq 0$, and $\a<0$).
\end{Df}

The functions which form ${\Lp}$ are known as the Laguerre
entire functions \cit{Il}.
Due to Laguerre and P\'olya (see e.g. \cit{Il}, \cit{Lev}),
we know that
\begin{Pn}
\label{4pn}
The family ${\cal L}^{+}$ consists of all entire functions which are
the polynomials possessing real nonpositive zeros only or their
uniform limits on compact subsets of $\
\C$.
\end{Pn}
We denote by ${\cal P}^{+}$ the set of polynomials belonging
to ${\cal L}^{+}$.
It should be pointed out that a function $f$, possessing the
representation (\ref{6}) with given $\a \in \R $,
belongs to ${\cal B}_{b}$ with $b>\vert \a \vert$, thus
it belongs to ${\cal A}_{\vert \a \vert }$ (but it may not belong to
${\cal B}_{\vert \a \vert }).$ It is also worth to remark that
every $f$ being of the form (\ref{6}) may
be written $f (z) = \exp (\a z) h(z) $, where $h$ is an entire
function of exponential type zero (for the details see
 e.g. \cit{BG1}, \cit{Lev}).

Consider the families
\be
\label{7}
{ \cal L}_{a} \ \sta \ {\cal L}\cap {\cal A}_{a}, \ \
{\cal L}^{\pm}_a \ \sta \ {\cal L}^{\pm }\cap {\AC}_a .
\ee
Obviously the latter definition of ${\cal L}_0 $ coincides with that
given by Definition \ref{1df}, and ${\cal L}_0 = {\cal L}_0^{+}$.

Now let us return to the example considered just after
Theorem \ref{n1tm}. The sequence $ \{z^n /n! \}$, as well as its
limit in $\EC$, $f\equiv 0$, belong to ${\LC}_0 $, thus
the implication $ {\rm (a)} \ {\rm and} \ {\rm (b)} \Rightarrow
{\rm (c)}$ fails also on ${\LC}_0 $. But this is the unique example
of $f\in \Lp$ of this type.
\begin{Tm}
\label{n2tm}
Every sequence of functions $\{f_n, n \in \N \}\subset
{\cal L}^+_{a}$, $a \geq 0$,
that converges in $\EC$ to a function $f\in {\AC}_a $,
which does
not vanish identically, is a bounded subset of this
${\AC}_a $ and hence converges to $f$ also in
${\AC}_a$. The limit function belongs to $\Lp$ as well.
\end{Tm}
\begin{Co}
\label{Lco}
A subset $B\subset {\cal L}^+_a $, which is compact in $\EC$,
is also compact in $\Aa$ provided it does not contain $f\equiv 0$.
\end{Co}

\begin{Co}
\label{n11co}
For every $a\geq 0$, the set ${\cal P}^+ $ is dense in
${\cal L}^{+}_{a}$ in ${\AC}_a $.
\end{Co}

Now we establish some sufficient conditions for a sequence in
${\LC}$ to be bounded in some ${\AC}$ not assuming its convergence.
Let $\{f_{n}, \ n\in \N \}$ be a sequence of functions from $%
{\cal L}.$ Denote by $C_{n},$ $l_{n},$ $\a_{n}$, and $\b_j (n)$
 the corresponding parameters of $ f_{n}$ in its representation
(\ref{6}). Let also
\begin{equation}
\label{8}
\mu_{k}(n)\ \sta \ \sum_{j=1}^{\infty }\b_{j}^{k}(n), \ \ k\in \N.
\end{equation}
\begin{Pn}
\label{5pn}
Given a sequence $\{ f_{n}, \ n\in \N \} \in {\cal L}$,
let there exist positive
$a$, $C$, and $l\in \N_0 $ such that, for all $n\in \N $,
$$
 \lv \a_n \rv +\mu_1 (n) \leq a ; \ \  \left| C_{n}\right| \leq C;
\  \  l_{n}\leq l.
$$
Then this sequence $\{f_n ,n \in \N \}$ is bounded
in ${\cal A}_{a}.$
\end{Pn}

{\bf Proof.} For $k< l_n$, $ f_n^{(k)} (0) = 0 $.
 For $k\geq l_n$, a simple calculation based
on the representation (\ref{6})
yields
\beq
f_{n}^{(k)}(0) & = &  C_{n} \frac{k!}{(k-l_{n})!}%
\sum_{i=0}^{k-l_{n}}\left[{ {k-l_{n}}\choose{i}}
 \a_{n}^{k-l_{n}-i} \right. \nonumber \\
& & \left. i!\sum_{1\leq j_{1}<j_{2}<...<j_{i}<\infty }
\b_{j_{1}}(n)\b_{j_{2}}(n)...\b_{j_{i}}(n)\right].  \nonumber
\eeq
Then
\begin{eqnarray}
\label{11}
\left| f_{n}^{(k)}(0)\right|
&\leq & Ck^{l_{n}}\sum_{i=0}^{k-l_{n}}{ {k-l_{n}}\choose{i}}%
\lv \a_{n}\rv^{k-l_{n}-i}(\mu_{1}(n))^{i} \nonumber \\
& = & C k^{l_n} (\lv \a_{n} \rv +\mu_{1}(n))^{k-l_{n}}\leq
C(k/a)^{l_{n}}a^{k}. \nonumber
\end{eqnarray}
 Thus for $b>a,$%
\begin{equation}
\label{12}
\left\| f_{n}\right\| _{b}=\sup_{k\in \N_{0}}\left\{ b^{-k}\left|
f_{n}^{(k)}(0)\right| \right\} \leq \sup_{k\in \N_{0}}\left\{
C(k/a)^{l_{n}}(a/b)^{k}\right\}.
\end{equation}
Finally, in view of the fact $l_{n}\leq l,$ the above estimate implies
that there exists
a constant $K$ depending only on $a$ and $b$ such that $\left\| f_{n}\right\|
_{b}\leq K.$
\kasten

\subsection{Operators}
For $\theta \geq 0$, let us consider a map $\Delta_{\theta}: %
\EC \ra \EC $ defined as follows
\be
\label{a1}
(\D_\th f)(z)  =  (\th + zD)Df(z) = \th \frac{df}{dz} + z \frac{d^2f}{dz^2}.
\ee
The study of
this map is quite important in view of the following facts. First,
one observes that, for $g (z) = f(z^2 )$,
\be
\label{aa01}
(\D_\th f)(z^2) = \frac{1}{4}\left(\frac{2\theta -1}{z}\frac{dg(z)}{dz}
+ \frac{d^2g (z)}{dz^2}\right),
\ee
which means that, for $\theta = N/2$, $N\in \N$, the map (\ref{a1})
is connected by the latter identity with the radial part
$\Delta_r$ of
the $N$--dimensional
Laplacean
$$
\D_r = \frac{N-1}{r}\frac{\pp}{\pp r}+ \frac{\pp^2}{\pp r^2}.
$$
This connection will be used by us in a separate work. Another
application arises from the fact that this map may produce
the Laguerre polynomials, which usually are defined as follows
(see e.g. \cit{Er}, p. 147)
\be
\label{LP1}
\tilde{L}_n^{(\th -1)} (z) \ \sta \ (-1)^n z^{-\th +1}e^z \left(
D^n z^{\th + n -1}e^{-z}\right),
\ee
namely
\beq
\label{LP2}
 \tilde{L}_n^{(\th -1)} (z)  =  e^z \Delta_\theta^n e^{-z}
                =  \exp( -\Delta_\theta)z^n .
\eeq
The latter formula contains an expression which in general
situations needs to be defined more precisely.

For given two entire functions $\varphi $ and $f$, we
denote by $\varphi (\Delta _{\theta })f(z)$ the formal series
\begin{equation}
\label{a2}
\sum_{k=0}^{\infty }\sum_{m=0}^{\infty }\frac{\varphi ^{(k)}(0)}{k!}\frac{%
f^{(m)}(0)}{m!}\Delta _{\theta }^{k}z^{m}.
\end{equation}
One can verify that
$$
\Delta _{\theta }^{k}z^{m} =q_{\theta }^{(m,k)}z^{m-k},\qquad q_{\theta
}^{(m,k)}=\left\{
\begin{array}{ll}
0, & k>m \\
\\
\gamma _{\theta }(m)/\gamma _{\theta }(m-k),  & 0\leq k\leq m
\end{array}
\right. ,
$$
where
\be
\label{a3}
\gamma_{\theta }(m)=m!\Gamma (\th + m ).
\ee

\begin{Pn}
\label{a1pn}
For positive $a$ and $b$ obeying the condition $ab<1$,
let $\vp\in {\BC}_a $ and $f\in {\BC}_b $. Then, for
every $\theta \geq 0$, the function
$g(z)=\varphi (\Delta _{\theta })f(z)$
belongs to ${\cal B}_{c}$ with $c=b(1-ab)^{-1}.$ Furthermore
\be
\label{a4}
\Vert g\Vert _{c}\leq (1-ab)^{-\theta }\Vert \varphi \Vert _{a}\Vert f\Vert
_{b}.
\ee
\end{Pn}

{\bf Proof} According to (\ref{a2})
$$
g(z)=\sum_{n=0}^{\infty }\frac{g^{(n)}(0)}{n!}z^{n},
$$
with
\be
\label{ab4}
g^{(n)}(0)=\sum_{k=0}^{\infty }\frac{n!}{k!(n+k)!}\varphi
^{(k)}(0) f^{(n+k)}(0) q_{\theta }^{(n+k,k)}.
\ee
But $\left| \varphi ^{(k)}(0) \right| \leq a^{k}\Vert \varphi \Vert _{a}$
and $%
\left| f^{(m)}(0) \right| \leq b^{m}\Vert f\Vert _{b}$ (see (\ref{1})).
For positive $a$ and $b$ obeying the condition $ab<1$, one may show that
$$
\sum_{k=0}^{\infty }\frac{n!}{k!(n+k)!}(ab)^{k}q_{\theta
}^{(n+k,k)}=(1-ab)^{-n-\theta },
$$
which yields in (\ref{ab4})
\beq
\left| g^{(n)} (0) \right| & \leq & \sum_{k=0}^{\infty }\frac{n!}{k!(n+k)!}\left|
\varphi ^{(k)}(0) \right| \left| f^{(n+k)}(0)\right| q_{\theta }^{(n+k,k)}
\nonumber \\
& \leq &
(1-ab)^{-\theta }\Vert \varphi \Vert _{a}\Vert f\Vert _{b}
\left(\frac{b}{1-ab}\right)^{n}. \nonumber
\eeq
Hence $g\in {\cal B}_{c}$ and the estimate (\ref{a4}) holds.
\kasten
\begin{Co}
\label{1co}
For all $\theta \geq 0,$ $a\geq 0$, and $b \geq 0$, such that $ab <1$,
$(\varphi ,f)\mapsto \varphi (\Delta_{\theta })f$ is a
continuous bilinear map from ${\cal A}_{a} \times {\cal A}_{b}$ into
${\cal A}_{c},$ where $c=b(1-ab)^{-1}$.
\end{Co}
The following strengthening
of the above statement is one of the main results of this research.
\begin{Tm}
\label{a1tm}
For all $\theta \geq 0,$ $a\geq 0$, and $b \geq 0$, such that $ab <1$,
 $(\varphi ,f)\mapsto \varphi (\Delta_{\theta })f$ is a
continuous map from ${\cal L}_{a}^{+}\times {\cal L}_{b}^{+} $ into
${\cal L}_{c}^{+},$ where $c=b(1-ab)^{-1}$.
\end{Tm}
An important kind of such operators corresponds to the choice of $\vp$
being of the form $\varphi_{a}(z)=\exp (az)$.
Every $\vp_a (\Delta_\th )$, $a\in \C$ maps
$\AC_0 $ into itself (see Corollary \ref{1co}).
Moreover, the family $\{ \vp_a (\Delta_\th )  \ ,
 a \in \ \C \}$, defined on $\AC_0 $, has a group property
\be
\label{ad21}
\vp_a (\Delta_\th )\vp_{a'} (\Delta_\th ) = \vp_{a+a'} (\Delta_\th ),
\ee
that may be proved on the base of (\ref{a2}). For
$a\geq 0$, the function $\vp_a $ belongs to $\Lp $.
In this case the operator $%
\varphi_a (\Delta_{\theta })$ has the following integral representation.
\begin{Pn}
\label{a2pn}
For every $\theta \geq 0,$ $a>0$, $b\geq 0$,
such that $ab <1$, and for all $f\in {\cal A}_{b}$,
\begin{eqnarray}
\label{a21}
(\exp(a \Delta_{\theta })f)(z) & = & \exp (-\frac{z}{a})
\int_{0}^{\infty}s^{\theta -1}e^{-s}w_{\theta }
\left(\frac{sz}{a}\right)f(as)ds \\
& = & \int_{0}^{\infty}K_\th \left(\frac{z}{a}s\right)f(as)
s^{\theta -1}e^{-s} ds \nonumber ,
\end{eqnarray}
with
\begin{equation}
\label{a22}
K_\th (z, s) \ \sta \ e^{-z}w_\th (zs), \ \
w_{\theta }(z)\ \sta \ \sum_{k=0}^{\infty }\frac{z^{k}}{\gamma _{\theta }(k)}.
\end{equation}
\end{Pn}
\begin{Rk}
\label{kernrk}
The integral kernel just appeared has the following expansion
in terms of the Laguerre polynomials (\ref{LP1}), (\ref{LP2})
$$
K_\th (z, s) = \sum_{n=0}^{\infty}\frac{z^n}{n!}
\frac{\tilde{L}_n^{(\th -1)}(s)}{\Gamma (\theta + n)}.
$$
Therefore,
$\{\tilde{L}_n^{(\th -1)}(s) / \Gamma (\theta + n) \}$
are the generalized Appell polynomials
with respect to the kernel $K_\theta $ (see \cit{BB}, p. 17).
\end{Rk}

{\bf Proof of Proposition \ref{a2pn}}.
 By Corollary \ref{1co}, $\exp(a\Delta_\th )$ is a
continuous operator on $\AC_b $, thus the left hand side
of (\ref{a21}) is well defined giving a function from
$\AC_c $, $ c= b(1-ab)^{-1}$. Then statement (i) of
Theorem \ref{n1tm} and the continuity of
the operator imply that the representation (\ref{a21})
needs  to be
proved only for $f(z)=z^{m}$. The definition (\ref{a2}) yields
\begin{equation}
\label{a23}
\exp (a{\Delta }_{\theta})z^{m}=
\sum_{n=0}^{m}\frac{a^{n}}{n!}\frac{\gamma _{\theta }(m)}{\gamma
_{\theta }(m-n)}z^{m-n}.
\end{equation}

It is not difficult to show that, for this choice of $f$, the
summation and the integration
in the right hand side of
(\ref{a21}) may be interchanged, which gives
$$
{\rm RHS}(\ref{a21})=\exp (-\frac{z}{a})\sum_{k=0}^{\infty }\frac{a^{m-k}z^{k}}{%
\gamma _{\theta }(k)}\Gamma (\theta +m+k)
$$
By means of the following Vandermonde--like
convolution identity (the proof see below)
\be
\label{aa23}
\frac{\Gamma (z+m+k)}{\Gamma (z+m)\Gamma (z+k)}=\sum_{n=0}^{\min (m,k)}%
{ {m}\choose {n}} {{k} \choose {n}}\frac{n!}{\Gamma (z+n)};\ z\in \ \C ; \ m,k,\in
\N,
\ee
we have:
\beq
{\rm RHS}(\ref{a21}) & = &\exp (-\frac{z}{a})\sum_{k=0}^{\infty }\frac{a^{m-k}z^{k}}{%
k!}\Gamma (\theta +m)\frac{\Gamma (\theta +m+k)}{\Gamma (\theta +m)\Gamma
(\theta +k)} \nonumber \\
& = & \exp (-\frac{z}{a})\sum_{k=0}^{\infty }\frac{a^{m-k}z^{k}}{k!}\Gamma
(\theta +m) \nonumber \\
&  &  \sum_{n=0}^{\min (m,k)}\frac{m!}{n!(m-n)!}\frac{k!}{n!(k-n)!}%
\frac{n!}{\Gamma (\theta +n)} \nonumber \\
& = & \exp (-\frac{z}{a})\sum_{n=0}^{m}\frac{a^{m-n}z^{n}m!}{n!(m-n)!}\frac{%
\Gamma (\theta +m)}{\Gamma (\theta +n)}\sum_{k=n}^{\infty }\frac{%
a^{n-k}z^{k-n}}{(k-n)!} \nonumber \\
& = & \sum_{n=0}^{m}\frac{a^{m-n}z^{n}m!}{(m-n)!n!}\frac{\Gamma (\theta +m)}{%
\Gamma (\theta +n)}= {\rm RHS}(\ref{a23})= {\rm LHS}(\ref{a21}) . \nonumber
\eeq
\kasten
The assertion just proved may be used to extend the
operator $\exp(a\Delta_\th )$.
In this case one ought to consider
the representation (\ref{a21}) as a
definition of the extended operator.
Here its following property -- a kind of the operation rule
(c.f. \cit{DTC}) --
 may be useful.
\begin{Pn}
\label{a3pn}
Given $a>0$ and $u\in \R$, let $b$ satisfy
$ 0\leq b < -u + 1/a$. Then, for every $g\in \AC_b $, the operator
(\ref{a21}) may be applied to the function
\be
\label{a24}
f(z) = \exp(uz)g(z),
\ee
yielding
\beq
\label{a25}
(\exp(a \Delta_\th )f)(z) & = & (1- ua)^{- \th }\exp\left(\frac{uz}{1-ua}\right)h(z),
\eeq
where
\beq
\label{ad25}
h(z)& = & \left[\exp\left(\frac{a}{1-ua}\Delta_\th \right) g \right] \left(\frac{z}{(1-ua)^2 }%
\right)  \\
& = & \exp[a(1-ua)\Delta_\th ]\left[g\left(\frac{z}{(1-ua)^2}\right) \right] . \nonumber
\eeq
Moreover, $h\in \AC_c $ with $c=b(1-ua)^{-1}[1-a(u+b)]^{-1}$.
\end{Pn}
\begin{Rk}
\label{a2rk}
For a negative $u$, the above statement extends the considered operator on
$\AC_d $ with $d < \vert u \vert + 1/a $, but this obviously
does not exhaust
all possible extensions.
The right hand side of (\ref{a21}) may be used to
define an integral operator
in the Hilbert space $L^2 (\R_{+}, \mu_\th )$
 possessing the kernel $K_\th (z, s)$
(\ref{a22}).
Here $\R_{+} = [0, +\infty) $ and $\mu_\th $, $\th >0 $ is the Euler measure
$$
\mu_\th (ds ) \ \sta \ \frac{1}{\Gamma(\th )}s^{\th -1 }e^{-s}ds .
$$
We construct such and other similar extensions in a separate work.
\end{Rk}
{\bf Proof of Proposition \ref{a3pn}.} For $f$ given by (\ref{a24}), one has
\beq
\label{a26}
(\exp(a\Delta_\th )f)(z)& & \\
  = &  & \exp(- \frac{z}{a}) \int_{0}^{+\infty}
s^{\theta -1}\exp[-s(1-ua)]w_{\theta }
\left(\frac{sz}{a}\right)g(as)ds. \nonumber
\eeq
A simple calculation yields
\beq
{\rm RHS} (\ref{a26}) & = &(1-ua )^{- \th} \exp(- \frac{z}{a}) \nonumber \\
& & \int_{0}^{+\infty}
s^{\theta -1}e^{-s}w_{\theta }\left(\frac{s(1-ua)}{a}
\frac{z}{(1-ua)^2}\right)
g(\frac{as}{1-ua})ds,
\nonumber \\
 & = &(1-ua )^{- \th} \exp(- \frac{z}{a})
\exp\left(\frac{1-ua}{a}\frac{z}{(1-ua)^2} \right) \nonumber \\
& &\left[\exp\left(\frac{a}{1-ua}\Delta_\th \right) g \right]
 \left(\frac{z}{(1-ua)^2 }\right),
\nonumber
\eeq
which gives
(\ref{a25}) and the first part of (\ref{ad25}). The second part of the latter
may be obtained by a change of variables. The final part of the statement follows
directly from Corollary \ref{1co}.
\kasten
Employing the extension by (\ref{a21}) we obtain
an extended form of Theorem \ref{a1tm}.
\begin{Tm}
\label{a2tm}
For every $a>0 $, $\th \geq 0$, the operator $\exp(a\Delta_\th )$
(\ref{a21}) maps:
\vskip.1cm
\noindent
\begin{tabular}{ll}
{\rm (i)}   & $\Lp_b $ into $\Lp_c $, with $0\leq b<1/a $
        and $c=b(1-ab)^{-1}$; \\
{\rm (ii)} & $\LM $ into $\LM$.
\end{tabular}
\end{Tm}
\begin{Rk}
\label{rka2tmrk }
The operator (\ref{a21}) acts on the whole $\LM$ with no
growth restrictions.
\end{Rk}
{\bf Proof of Theorem \ref{a2tm}.}
Claim (i) is simply a repetition of Theorem \ref{a1tm},
which we add here in order to describe the action of this
operator on the whole $\LC$ in one statement.
 To prove (ii) we observe that every $f\in \LM$ may be
written in the form (\ref{a24}) with $g\in \LC_0 $ and $u<0$. Then, for every
positive $a$, one may choose $b =0$ and apply the operator (\ref{a21}) in accordance
with Proposition \ref{a3pn}. The result will be given by (\ref{a25}). Since
$a(1-ua)^{-1} $ is positive, $u(1-ua)^{-1}$ is negative and the function
$$
\left[\exp\left(\frac{a}{1-ua}\Delta_\th \right) g \right] \left(\frac{z}{(1-ua)^2 }\right)
$$
belongs to $\LC_0$, which means that its product with $\exp(u(1-ua)^{-1}z)$
belongs to $\LM$.
\kasten

Now we may use the operators introduced above
to solve the following initial value
problem.
\beq
\label{e1}
\frac{\partial f (t,z)}{\partial t} & = &\th \frac{\partial f(t,z)}{\partial z} +
z\frac{\partial^2 f(t, z)}{\partial z^2}, \ \ \ t\in \R_{+} , \ z\in \ \C , \\
f(0, z) & = & g(z) . \nonumber
\eeq
\begin{Tm}
\label{e1tm}
For every $\th \geq 0$ and  $g\in \EC $ having the form
\be
\label{e3}
g(z) = \exp(-\ve z)h(z), \ \ \ h\in {\AC}_0, \ \ \ve\geq 0,
\ee
\vskip.1cm
\begin{tabular}{ll}
{\rm (i)}  & the problem (\ref{e1}) has in ${\AC}_\ve $ the
following solution
\end{tabular}
\vskip.1cm
\beq
\label{e2}
f(t,z) & = & \left(\exp(t\Delta_\th )g \right) (z) \\
   & = & \exp(-\frac{z}{t})\int_{0}^{+\infty }
s^{\th -1} w_{\th}\left(\frac{z s}{t}\right)
 e^{-s}g(t s)d s, \ \ t>0 .
\nonumber
\eeq
\vskip.1cm
\begin{tabular}{ll}
     {\rm (ii)}
& If the initial condition $g$ possesses
 (\ref{e3}) with $\ve >0$,\\
  &  then the solution (\ref{e2}) converges to zero
  when $t \ra +\infty$ \\
& uniformly on compact
 subsets of $\ \C$.\\ [.1cm]
 {\rm (iii)} &If in (\ref{e3}) $h \in \Ln\subset{\AC}_0$,
 then the solution (\ref{e2}) \\
& belongs either to $\Ln$, for $\ve =0$,  or to $\LM$, for $\ve >0$.
\end{tabular}
\end{Tm}

Claim (ii) means that the so called stabilization of the solutions
holds (see e.g. \cit{Kam} and \cit{De}).

\section{Proofs}

\subsection{Spaces of entire functions }

We start with the proof of Theorem \ref{n1tm} by introducing
another norms.
For appropriate $f\in \EC$ and some $b>0$, we set
\be
\label{n1}
N_b (f) \ \sta \ \sup_{z\in \ \C}\{\lv f(z) \rv \exp(-b\lv z \rv ) \}
= \sup_{r\in \R_{+} } \{M_f (r)\exp (-b r) \},
\ee
where
$$
M_f (r) \ \sta \ \sup_{\lv z \rv \leq r } \lv f(z) \rv , \ \
r\in \R_{+} .
$$
Obviously $N_b (.) $ is a norm on a subset of $\EC$.
\begin{Pn}
\label{n1pn}
For given $f\in {\BC}_b $, let $N_{b-\ve} (f) < \infty $ with some
$\ve \in (0, b)$. Then there exists a constant
$C(b , \ve )$ such that
\be
\label{n2}
N_b (f) \leq \|f \|_b \leq C(b,\ve )N_{b-\ve }(f).
\ee
\end{Pn}
{\bf Proof.}
By means of the Cauchy inequality one obtains
$$
\lv f^{(k)} (0) \rv \leq \frac{k!}{r^k }M_f (r) , \ \ k\in \N , \ r\in \R_{+}.
$$
By the definition
$$
M_f (r) \leq N_{b-\ve }(f)\exp[(b-\ve)r] ,
$$
thus
\be
\label{n3}
\lv f^{(k)} (0) \rv b^{-k} \leq  N_{b-\ve }(f) k! \chi_k (r) .
\ee
The function $\chi_k (r) \ \sta \ (b r)^{-k}\exp[(b-\ve)r]$ has
the unique minimum at $r=k(b-\ve)^{-1}$, hence the latter estimate
would be the best possible for this value of $r$. We set
$$
C_0 (b , \ve) =1 , \ \ C_k (b, \ve ) = k!
\chi_k \left(\frac{k}{b-\ve}\right) =
\frac{k!}{k^k} \left(1-\frac{\ve}{b}\right)^k e^k  ,
$$
and obtain in (\ref{n3})
$$
\lv f^{(k)} (0) \rv b^{-k} \leq C_k (b, \ve ) N_{b-\ve }(f) .
$$
By means of the Stirling formula, one may get convinced that the
sequence $\{C_k (b , \ve )\}$ is bounded.
Thus we set
$$
C(b , \ve) \ \sta \ \sup_{k\in \N_0 }C_k (b, \ve ) ,
$$
and obtain the upper bound of $\|f\|_b $ in (\ref{n2}).
To complete the
proof we observe that
$$
M_f (r) \leq \sum_{k=0}^{\infty}
 \frac{1}{k!}
\lv f^{(k)} (0) \rv b^{-k}
(b r)^k \leq \|f\|_b \exp(b r) ,
$$
which immediately yields the lower bound in (\ref{n2}).
\kasten

The family
$\{N_{b} (.) \ b > a \}$ defines a topology on
${\AC}_a $, which by Proposition \ref{n1pn} is equivalent to
the topology ${\cal T}_a $ introduced above.
 We use this fact as follows.

\begin{Pn}
\label{n2pn}
For every $a\geq 0 $, the space ${\cal A}_{a}$ possesses the
properties:
\vskip.1cm
\begin{tabular}{ll}
{\rm (i)}  & let $f\in {\cal A}_{a}$ and $g\in {\cal A}_{b}$, then
            their product $f g$ belongs to ${\cal A}_{a+b}$; \\ [.1cm]
{\rm (ii)} & let $f \in {\cal A}_{a}$ and, for some $g\in \EC$,
            the function \\
           & $\l (r, r_0 )\ \sta \ \log M_g (r) - \log M_f (r +r_0 ), $
             with some
             fixed \\
           & $r_0 \in \R_{+}$, be bounded as a function of $r\in \R_{+}$,   \\
            & then $g$ also belongs to ${\cal A}_{a}$.
\end{tabular}
\end{Pn}
{\bf Proof.}
Here we define the topology ${\cal T}_a $
by means of the
family $\{ N_{b} (.) , \ b > a \}$.  The proof of (i) is obvious.
The proof of (ii) is also quite simple:
\beq
N_{b} (g) & = & \sup_{r\in \R_{+}} \{ M_f (r+ r_0 )
\exp [ - b(r+r_0 )]
\exp[ b r_0 + \l(r, r_0 )] \} \nonumber \\
&\leq & N_{b} (f) \exp( b r_0 ) \sup_{r \in \R_{+}}
 \exp[\l (r, r_0 )] . \nonumber
\eeq
\kasten
{\bf Proof of Theorem \ref{n1tm} and Corollaries \ref{pMonco} --
\ref{new1co}.}
Proposition \ref{n2pn} yields that we may use here
 Proposition 2.7 and Proposition 2.5 of \cit{Tay}, which imply
claim (i) and claim (ii) respectively. Corollary \ref{pMonco}
follows directly from claim (ii), which also makes possible
to extend Montel's property on $\Aa$ from $\EC$.
The proof of Corollary \ref{new1co} may be given as follows.
Both sequences are bounded in the corresponding spaces.
By means of the triangle inequality one gets
$$
 N_{a' + b'}(f g - f_n g_n ) \leq
N_{a'}(f_n )N_{b'}(g_n - g ) + N_{a'}(f_n - f)N_{b'}(g ),
$$
which yields the converges to be proved.

\kasten

Now let $\{f_n , \ n\in \N \} \subset\Lp$
converge in $\EC$ to a function $f$,
which does not vanish identically.
Then by Proposition \ref{4pn} the latter function
also belongs to $\Lp$ and each such a function may be written
in the form (\ref{6}). Since some of the negative zeros
of $f_n$ could converge
to zero, certain sequences $\{\b_j (n), \ n\in \N\}$
would be unbounded and at the same time the sequence $\{C_n \}$
would converge to zero.
In view of this possibility
it is more convenient to rewrite (\ref{6})
for $f_n $ as follows
\be
\label{ne4c}
f_n (z) =   p_n (z)\tilde{f}_n (z), \ \
p_n (z)  =  C_n z^{l_n }\prod_{k=1}^{q_n}(z+z_k (n)) ,
\ee
\be
\label{ne4d}
 \tilde{f}_n (z)  =  \exp(\a_n z)\prod_{j=1}^{\infty}
(1+ \b_j (n) z) , \nonumber \\
\ee
and suppose that the sequences
 $\{\b_j (n) ,n\in\N \} $ are bounded and the sequences
$\{z_k (n) ,n\in\N \}$ converge to zero. We also
write
\beq
\label{ne4a}
f(z)& = &  p(z)\tilde{f}(z) , \ \ \ p(z)  =  Cz^l , \\
\label{ne4b}
  \tilde{f}(z) & = & \exp(\a z)\prod_{j=1}^{\infty}
(1+ \b_j  z) .
\eeq
Recall that in the above representation all $\b_j (n)$ and
$\b_j $ are numbered according to the definition (\ref{6}),
i.e., $\b_j (n) \leq \b_{j+1}(n)$
The following statement, which describe the convergence
of the sequences $\{C_n \}$,  $\{\b_j (n)\}$,
$\{z_k (n)\}$, follows directly from the assumed convergence
$f_n \ra f$ by known Hurwitz's theorem (see \cit{BG}, p. 167).
\begin{Pn}
\label{n24pn}
There exist positive integers $n_{*}$, $l_{*}$, and $ q_{*}$
such that $l_{*} + q_{*} = l$ and, for all $n>n_{*} $, $l_n = l_{*}$,
$q_n = q_{*}$,
all the sequences $\{z_k (n), \ n= n_{*}+1, \dots \} $, $k = 1, \dots
q_{*}$ converge to zero, and the sequence
$\{\b_1 (n), \ n= n_{*}+1, \dots \} $ converges to $\b_1 $.
If all $\b_j = 0$, then all the sequences
$\{\b_j (n) , \ n\in \N \}$ converge to zero.
\end{Pn}
{\bf Proof of Theorem \ref{n2tm}.}
By Proposition \ref{4pn}
the convergence $f_n \ra f$ in $\EC$ implies that $f\in \Lp$,
which proves the final part of the theorem.
For the functions considered, we use the forms (\ref
{ne4c}) -- (\ref{ne4b}). The convergences described
by Proposition \ref{n24pn} implies
$p_n \ra p $ in $\EC$. As a sequence
of polynomials of bounded degree,
the sequence $\{p_n \}$ is bounded
in ${\AC}_0 $ hence it converges to $p$ in ${\AC}_0 $.
Therefore, by Corollary \ref{new1co} it remains to prove
the following convergence in ${\AC}_a $
\be
\label{ne5}
\tilde{f}_n   \ra \ \tilde{f} .
\ee
As for the parameter $\a $ in (\ref{ne4a}), it must be bounded
$\a \leq a$ since $f \in {\AC}_a$.
 We choose $\a = a$
-- the convergence (\ref{ne5}) for the
other values of $\a$ will follow from the proof for this choice.
The above proven convergence $p_n \ra p$ and the assumed convergence
in $\EC$ of the sequence $\{f_n \}$ imply also
the convergence $\tilde{f}_n   \ra \ \tilde{f}$ in $\EC$, which
yields by Weierstrass' theorem
\be
\label{k1}
\tilde{f}_n^{(k)} (0) \ra \tilde{f}^{(k)} (0), \ \ k\in\N_0 ,
\ \ n\ra \infty.
\ee
For all $n\in\N$,
$\tilde{f}_n (0) = \tilde{f}(0) =1$, therefore,
the functions $\psi_n (z) \ \sta \ \log \tilde{f}_n (z)$,
$\psi (z) \ \sta \ \log \tilde{f} (z)$
are differentiable at $z=0$,
their derivatives of order $k \in \N$ may be written as
polynomials of corresponding $f_n^{(l)} (0)$ or $f^{(l)} (0)$
with $l = 0,1, \dots , k$. Hence (\ref{k1}) yields
$$
\psi_n^{(k)} (0) \ra \psi^{(k)} (0), \ \ k\in \N, \ \ n\ra \infty,
$$
which may be written
\be
\label{k2}
\mu_k (n) \ra \mu_k \ \sta \ \sum_{j=1}^{\infty}\b_j^k , \ \ k\geq 2 ,
\ee
and
\be
\label{k3}
  \a_n + \mu_1 (n) \ra a+ \mu_1  ,
\ee
where $\mu_1 $ is defined by (\ref{k2}) with $k=1$ and
the notation (\ref{8}) has been used.
In view of claim (ii) of Theorem \ref{n1tm} we will show the
convergence (\ref{ne5}) by proving the boundedness of the
sequence $\{\tilde{f}_n \}$ in ${\AC}_a$. To this end we use
the topology
on ${\AC}_a$ defined by the family of norms (\ref{n1}).
It is seen that, for $c>a$,
\be
\label{ne6}
N_c (\tilde{f}) =  \sup_{r\in \R_+}\{\exp[-(c-a)r]\prod_{j=1}^{\infty}
(1+ \b_j  r)\}.
\ee
Since all $f_n \in \Aa$, the sequence of nonnegative parameters
$\{\a_n \}$ is bounded $\a_n \leq a$.
Similarly to (\ref{ne6}), we obtain
\be
\label{ne7}
N_c (\tilde{f}_n) =  \sup_{r\in \R_+}\{\exp[-(c-\a_n )r ]
\prod_{j=1}^{\infty}
(1+ \b_j (n) r)\}, \ \ c>\a_n.
\ee
First we consider the simplest situations. Suppose that all
$\b_j (n)$ and $\b_j $ equal to zero. Then (\ref{k3}) implies
$\a_n \ra \a$ and (\ref{ne5}) obviously holds.
Suppose now that all $\b_j =0$ and all $\tilde{f}_n $,
except maybe a finite number of such functions,  have
finitely many $\b_j (n)$ different from zero. Then
by means of Proposition \ref{n24pn} one can easily
show that
$$
\prod_{j=1}^{m_n}(1+ \b_j (n) z ) \ra 1, \ \ \ n\ra \infty
$$
in ${\AC}_0 $,
as it took place with the convergence $p_n \ra p$.
Then again $\exp( \a_n z ) \ra  \exp( a z )$ in $\Aa$.
The remaining situations are more complicated.

The boundedness of the sequence $\{\tilde{f}_n\} $ in $\Aa$
is proven by showing that the sequence of norms
$\{\N_c (\tilde{f}_n )\}$ is bounded for all $c>a$.
Let us start with the evaluation of
$N_c (\tilde{f})$. If all $\b_j =0 $, then $N_c (\tilde{f}) =1$
for all $c>a$.
 In the case of nonzero $\b_j$, we
find the point $r_c \in \R_+ $ where the supremum in (\ref{ne6})
is achieved. It may be done by solving
the equation
\be
\label{ne8}
c = a + \sum_{j=1}^{\infty}\frac{\b_j}{1+ \b_j r} \ \sta \ \vp (r).
\ee
Except for the case where all $\b_j =0 $, which we consider
at the end of this proof,
 $\vp $ is a monotone decreasing function
on $\R_+ $. Since the series $\sum \b_j $ converges (see (\ref{6})),
the second term in
(\ref{ne8}) tends to zero when $r\ra +\infty$, hence $\vp $ takes
 on $\R_+ $ all values from $(a, a+\mu_1 ]$.
Thus, for $c>a+\mu_1 $,
$N_c (\tilde{f}) = 1$. For $c \in (a, a+\mu_1]$,
\be
\label{ne9}
 N_c (\tilde{f}) = \exp[-(c-a)r_c]\prod_{j=1}^{\infty}
(1+ \b_j  r_c),
\ee
where $r_c $ is the unique solution of the equation (\ref{ne8}).
Similarly
one obtains from (\ref{ne7}) the following equation
\be
\label{ne10}
c = \a_n + \sum_{j=1}^{\infty}\frac{\b_j (n)}{1+ \b_j (n) r}\
 \sta \ \vp_n (r),
\ee
where $\vp_n $ is also a monotone decreasing function taking on
$\R_+ $ all values from the interval $(\a_n, \a_n + \mu_1 (n)]$.
Thus, for $c>\a_n + \mu_1 (n) $, $ N_c (\tilde{f}_n) =1$. For
$c\in (\a_n, \a_n + \mu_1 (n)]$, the equation
(\ref{ne10})
has the unique solution $r_n$, which defines the norm
\be
\label{ne11}
 N_c (\tilde{f}_n) = \exp[-(c-\a_n ) r_n ]\prod_{j=1}^{\infty}
(1+ \b_j (n) r_n) .
\ee

Obviously each $\vp_n $ may be analytically
continued on the complex half-plane
$A_\ve \ \sta \ {\rm Re}r \geq -\ve $ with some
$\ve >0 $, which obeys the conditions
\be
\label{ve}
 \ve \sup_{n> n_* }\b_1 (n) \ \sta \ \ve \b   <1 .
\ee
Such supremum exists by Proposition
\ref{n24pn}. Let us show that the sequence $\{\vp_n , \ n=
n_* +1 , \dots \}$ is bounded on $A_\ve $. Set
$r+ \ve = z= x+iy$, then
\be
\label{ve1}
\vp_n (r) = \a_n + \sum_{j=1}^{\infty}
\frac{\b_j^* (n)}{1+ \b_j^* (n)z}, \ \  \b_j^* (n) \ \sta  \
\frac{\b_j (n)}{1- \ve \b_j (n)} .
\ee
Further, for $r\in A_\ve $, $x\geq 0$ and $y\in \R$
and  one readily obtains
\beq
 0 & < &{\rm Re}\vp_n (r)  =  \a_n + \sum_{j=1}^{\infty}
\frac{\b_j^* (n)(1+\b_j^* (n)x ) }{(1+ \b_j^* (n)x)^2 +
(\b_j^* (n)y)^2}  \nonumber \\
& \leq & \a_n + \sum_{j=1}^{\infty}
\frac{\b_j^* (n) }{1+ \b_j^* (n)x } \leq
\a_n + \sum_{j=1}^{\infty}\b_j^* (n) \nonumber \\
& \leq & a +
\frac{\mu_1 (n) }{1- \ve \b} \leq a +
\frac{1  }{1- \ve \b} \sup_{n> n_* }\mu_1 (n). \nonumber
\eeq
In view of (\ref{k3}) the latter supremum exists.
Similarly
\beq
\lv {\rm Im}\vp_n (r) \rv & = & \sum_{j=1}^{\infty}
\frac{[\b_j^* (n)]^2 \lv y \rv }{(1+ \b_j^* (n)x)^2 +
(\b_j^* (n)y)^2}
\leq \sum_{j=1}^{\infty}
\frac{\b_j^* (n) \lv y \rv }{1+
(\b_j^* (n)y)^2} \b_j^* (n) \nonumber \\
& \leq &  \sum_{j=1}^{\infty}\b_j^* (n) \leq
\frac{1  }{1- \ve \b} \sup_{n> n_* }\mu_1 (n).\nonumber
\eeq
Therefore, by Montel's compactness criterium (\cit{BG}, p.120)
the sequence $\{\vp_n \}$ is relatively compact on $A_\ve $.
It is easily seen that the convergence established
by (\ref{k1}), (\ref{k2}) implies $\vp_n^{(k)} (0)
\ra \vp^{(k)} (0)$ for all $k\in \N_0 $. This is enough
for the convergence
$\vp_{n }\ra \vp $ and also for the derivatives
$\vp_n^{\prime} \ra \vp' $,
uniformly on compact subsets of $A_\ve $ (see \cit{BG}, p.121).
But $\vp' $ is strictly negative for all $r\in \R_+ $, which
means that, for any $R>r_c $, there exists $C(R) $ such that
$$
0< C(R) \leq \inf_{r\in [0,R]}\lv \vp_n^{\prime} (r) \rv ,
$$
for all sufficiently large $n$. Here
we assume that $n_*$ is such that
this estimate holds for all $n>n_*$.
A simple use of the proven uniform convergences yields
that $r_n \ra r_c $ with the estimate
$$
\lv r_n - r_c \rv \leq \frac{1}{C(R)}\sup_{r\in [0,R]}
\lv \vp_n (r) - \vp (r) \rv .
$$
This yields the boundedness of the sequence of
norms $N_c (\tilde{f}_n )$, that
was to be proven. It remains to consider the case where
$\tilde{f} (r) = \exp(a r)$. Here $\vp (r) \equiv a$ and
by Proposition \ref{n24pn}
all $\b_j (n) $ tend to zero. Moreover, (\ref{k2}) and (\ref{k3})
yield in
this case
$\a_n + \mu_1 (n) \ra a$ and $\mu_k (n) \ra 0$, $k\geq2$.
Then we write
$$
\vp_n (r) = \a_n +  \mu_1 (n) + \om_n (r),
$$
with
$$
\left\vert \om_n (r) \right\vert \leq \max\{r\mu_2 (n),
r^2 \mu_3 (n)\},
$$
which yields $\vp_n (r) \ra a$ uniformly on every $[0,R]$. This
means $N_c (\tilde{f}_n ) \ra 1$.
\kasten
{\bf Proof of Corollaries \ref{Lco} , \ref{n11co}.}
A compact subset of $\EC$ is bounded and closed
in $\EC$.
Then, being a subset of $\Lpa$, it is bounded
and closed in $\Aa$ provided it does not contain $f\equiv 0$
by Theorem \ref{n2tm}.
Hence it is compact in $\Aa$ in view of its Montel's property.
By Proposition \ref{4pn}, every sequence of
polynomials from ${\cal P}^+ $ converges in $\EC$ to some
$f\in \Lpa $ $a\geq 0$, and every such a function is a limit in $\EC$
of a sequence of polynomials from ${\cal P}^+ $. The sequence
of polynomials
obviously is a subset of $\Aa$ with any $a\geq 0$. By the above
theorem this sequence converges to $f$ in $\Aa$.
\kasten

\subsection{Operators}

 The proof of Theorem \ref{a1tm} is divided on several steps.
Below we will need
a tool to control the distribution of zeros of certain holomorphic functions.
Thus we begin
with its construction. Introduce
\begin{equation}
\label{aa4}
A=\{ z\in \C \ \vert \ {\rm Re}z >0 \},\ \
 \ \bar{A}=\{z\in \C \ \vert \ {\rm Re}z \geq 0 \}.
\end{equation}
\begin{Lm}
\label{a1lm}
Let $Q_{0}$ and $Q_{1}$ be respectively a holomorphic
function in $\bar{A}$ and a polynomial in a single complex variable. If
\begin{equation}
\label{a5}
R(v,w) \ \sta \ Q_{0}(w)+vQ_{1}(w)\neq 0,
\end{equation}
whenever $v,$ $w\in \bar{A},$ then
\begin{equation}
\label{a6}
S(z) \ \sta \ Q_{0}(z)+{Q}_{1}^{\prime}(z)\neq 0,
\end{equation}
whenever $z\in \bar{A}.$
\end{Lm}

{\bf Proof. }
It is no need to prove the statement in the trivial case $Q_1 \equiv 0$. In the
nontrivial case the assumed property of the function $R$ yields that
$Q_{0}(w)= R(0, w) $ does not vanish whenever $w\in \bar{A}$. For
$v\in \bar{A}\setminus \{0\} $,
we rewrite (\ref{a5}) as follows
\be
\label{a7}
R(v,w)=v T(v^{-1},w),  \ \  T(\varepsilon ,w)=\varepsilon Q_{0}(w)+Q_{1}(w) .
\ee
Then $T(\varepsilon
,w)\neq 0$, for $\varepsilon \in \bar{A}\backslash \{0\}$ and $w\in \bar{A}.$
The above facts together with Rouch\'{e}'s theorem
(\cit{BG}, p. 167) imply that
$Q_{1}(w)\neq
0 $ for $w\in A.$ Let us decompose $\bar{A}$ onto the following
subsets
\be
\label{ab7}
B_{0}=\{z\in \bar{A} \ \mid \  Q_{1}(z)=0\},\ \  B_{1}=\bar{A}\setminus
B_{0}.
\ee
Recall that $Q_1 $ is a
polynomial, which does not vanish on $A$, thus $B_0 $ is a part
(maybe empty) of the finite set of $Q_1$ zeros, that is
\begin{equation}
\label{ac7}
B_{0} = \{ z_{1} ,
 \dots , z_{m }\}, %
 \ \ \  {\rm Re} z_{j } = 0 , \ %
j = 1, \dots , m .%
\end{equation}
For $ z\in \bar{A}$, we set
\be
\label{a8}
S(z)=Q_{s}(z)T_{s}(z),\ \ \ \ %
z \in B_{s}; \  s=0,\ 1;
\ee
with
$$
T_{0}(z)=1+\frac{Q_{1}^{\prime }(z)}{Q_{0}(z)},\ \ \ %
T_{1}(z)=\frac{Q_{0}(z)%
}{Q_{1}(z)}+\frac{Q_{1}^{\prime }(z)}{Q_{1}(z)}.%
$$
Let us prove now that $T_{s}(z)\neq 0,$ when $ z\in B_{s}.$

(a) For every  $z\in B_0$, there exist two possibilities:
\be
\label{ab8}
{\rm (i) } \ \ {\rm Im} \frac{Q_{1}^{\prime }(z)}{Q_{0}(z)} \neq 0 ; \ \
\ \ {\rm (ii) } \ \ {\rm Im} \frac{Q_{1}^{\prime }(z)}{Q_{0}(z)} = 0 .
\ee
The first one immediately yields $T_{0}(z)\neq 0$.
In the second case we show that the real part of $T_0 $ is not less
than 1. To this end
let us consider the following equation
\be
\label{ac8}
T(\varepsilon ,w)=0;\ \ \ \varepsilon \in \C , \ \ w\in \bar{A} ,
\ee
or
$$
\varepsilon Q_{0}(w)=-Q_{1}(w);\ \ \ \varepsilon \in \C , \ \ w\in
\bar{A}.
$$
Recall that $Q_0 $ is a holomorphic function on $\bar{A}$,
which does not vanish there. Hence, for $w\in {\bar A}$, the
latter equation has the following solution
$$
\ve = \ve (w) \ \sta \ - \frac{Q_1 (w)}{Q_0 (w) } ,
$$
which is a holomorphic function at any point of $\bar{A}$,
possessing isolated zeros in $ B_0 $.
Therefore, for every $z \in B_0 $, there exists a neighborhood
of this point where one may write
\beq
\varepsilon (w) & = & - \frac{Q_{1}^{\prime }(z)}{Q_{0}(z)}(w-z)+o(\left|
w-z\right| ). \nonumber
\eeq
Having in mind (see (\ref{ac7}))
that, for $z\in B_0 $, ${\rm Re} z = 0$ and that the second possibility
in (\ref{ab8}) is considered, one gets
\beq
\label{aa9}
{\rm Re}\ve (w) & = & - {\rm Re}\frac{Q_1^{\prime} (z) }{Q_0 (z) } {\rm Re}w
+ o(\vert w- z \vert ).
\eeq
From the assumption of this lemma and from the
 definition (\ref{a7}) we know that $T(\ve , w)$ does not vanish
whenever $\ve \in \bar{A}\setminus \{0 \}$ and $w\in \bar{A}$.
On the other hand, $\ve (w) $ is
a solution of the equation (\ref{ac8})
hence the values of the function $\ve (w)$ on $B_1$ should
have negative real parts only, i.e.
$$
{\rm Re} \ve (w) < 0 , \ \ \ {\rm for} \ \ w\in B_1  .
$$
The latter yields in turn in (\ref{aa9})
\begin{equation}
\label{a11}
0\geq \lim_{w\rightarrow z}\frac{{\rm Re}\varepsilon (w) }{{\rm Re}w}%
=-{\rm Re}\frac{Q_{1}^{\prime }(z)}{Q_{0}(z)} .
\ee
Thus
$$
{\rm Re}T_0 (z) = T_0 (z)=
1+ {\rm Re}\frac{Q_{1}^{\prime }(z)}{Q_{0}(z)} \geq 1,
$$
that was to be shown in the case (ii) in (\ref{ab8}).
In what follows, in both cases
$T_{0}(z)\neq 0$ whenever $z\in B_0$.

(b) Now we prove that $T_1 $ does not vanish on $B_1$. To this end we
write $T_{1}(z)=\vartheta (z)+t (z),$ with
$$
\vartheta (z)=\frac{Q_{0}(z)}{Q_{1}(z)};\ \ \ t (z)=\frac{Q_{1}^{\prime
}(z)}{%
Q_{1}(z)}=\sum_{j=1}^{M}\frac{1}{z- z_{j}}=\sum_{j=1}^{M}\frac{\bar{z}-%
\bar{z}_{j}}{\left| z-z_{j}\right| ^{2}} ,
$$
where $z_j $, $j = 1, \dots , M$ belong
either to $B_0$ (see  (\ref{ac7}))
or to $\ \C \setminus {\bar A} $. In any case
${\rm Re}z_{j}\leq 0,$ thus, for $z\in B_{1}$, $%
{\rm Re}t(z)\geq 0.$ Hence it suffices to prove that ${\rm Re}\vartheta(z)>0$
whenever $z\in B_{1}.$ Rewrite (\ref{a5}) in the form
\begin{equation}
\label{a12}
R(v,w)=Q_{1}(w)[\vartheta (w)+v];\ \ \ v\in \bar{A},\quad w\in B_{1}.
\end{equation}
Suppose ${\rm Re}\vartheta (w)\leq 0$, for some $w \in B_1 $. Then one may set
in (\ref{a12}) $v=-{\rm Re}\vartheta(w)-i{\rm Im}\vartheta (w)$
and obtain that
$R(v, w)=0,$ for $v\in {\bar A}$ and $w\in B_1 \subset {\bar A}$, which is
contradictory to the assumption
(\ref{a5}). That means $T_{1}(z)\neq 0.$
\kasten

\begin{Rk}
\label{a1rk}
The above lemma is a generalization of the similar statement proved by
E. H. Lieb and A.
D. Sokal in \cit{LiS}, where the case with both $Q_0$ and $Q_1 $ being
polynomials was considered. We use the possibility to
take $Q_0 $ being a meromorphic function below.
\end{Rk}
First we prove a simple corollary of Lemma \ref{a1lm}
\begin{Lm}
\label{a2lm}
Let $Q_{0}$ and $Q_{1}$ be respectively a holomorphic
function in $A$ and a polynomial in a single complex variable. If
\begin{equation}
\label{a13}
R(v,w)\ \sta \ Q_{0}(w)+vQ_{1}(w)\neq 0,
\end{equation}
whenever $v,$ $w\in A,$ then
\begin{equation}
\label{a14}
S(z)\ \sta \ Q_{0}(z)+Q_{1}^{\prime }(z)\neq 0,
\end{equation}
whenever $z\in A,$ or else $S(z)\equiv 0.$
\end{Lm}

{\bf Proof} For arbitrary $\delta >0$, we set
$$
R_{\delta }(v,w) \ \sta \ R(v+\delta ,w+\delta );\ \ \ S_{\delta  }(z) \ \sta
\ S(z+\delta).
$$
Clearly $R_{\delta }(v,w)\neq 0$ whenever $v,$ $w\in \bar{A},$ then
by Lemma \ref{a1lm} one gets $S_{\delta }(z)\neq 0, $ whenever $z\in \bar{A}%
. $ For $\delta \searrow 0,$ $S_{\delta }(z)$ is uniformly
convergent on compact subsets of $A$ to $S(z).$ Thus Hurwitz's theorem yields
$S(z)\neq 0$ on $A,$ or else $S(z)\equiv 0 $.
\kasten
Now we use the mentioned possibility to take $Q_0 $ being a meromorphic function on
$\bar{A}$.
\begin{Lm}
\label{a3lm}
Let $P$, $Q$, and $Q_{1}$ be polynomials in a single complex variable. Suppose
that $P $ does not vanish on $A$ and
\begin{equation}
\label{a15}
Q(w)+P(w)v Q_{1}(w)\neq 0,
\end{equation}
whenever $v,$ $w\in A.$ Then either
\begin{equation}
\label{a16}
S(z) \ \sta \ Q(z)+P(z)Q_{1}^{\prime }(z)\neq 0,
\end{equation}
whenever $z\in A,$ or else $S(z)\equiv 0.$
\end{Lm}

{\bf Proof } Since $P$ does not vanish on
$A$, we get
from (\ref{a15})
$$
\frac{Q(w)}{P(w)}+vQ_{1}(w)\neq 0 ,
$$
for $v$, $w$ belonging to $A$.
Setting $Q_{0}(w)=Q(w) /P(w) ,$ we get from Lemma \ref{a2lm}
that either
$$
\hat{S}(z)\ \sta \ Q_{0}(z)+Q_{1}^{\prime }(z)\neq 0,
$$
whenever $z\in A,$ or else $\hat{S}(z)\equiv 0.$ The
former implies (\ref{a16}).
\kasten

Now we are at a position to study our map $\D_{ \theta}$
\begin{Lm}
\label{a4lm}
For arbitrary nonnegative $\k$ and $\theta $,
$\k +\Delta _{\theta }$ maps $\cal{L}^{+}$ into
$\Lp$.
\end{Lm}

{\bf Proof.} For the continuous operator $\k
+\Delta _{\theta } : \EC \ra \EC$, it suffices to prove the stated property
on a ${\cal T}_C $--dense subset of $\LC^{+}$. A proper choice of such subset is
${\cal P}^{+}$. Then the statement of the lemma
is equivalent to the claim that the polynomial
$(\k +\Lambda _{\theta })q(z)$
with
\begin{equation}
\label{a17}
\Lambda _{\theta }=\left( \theta +\frac{z}{2}%
D\right)\left( \frac{1}{2z}D\right)  ,\ \ \ q(z)=p(z^{2}),
\end{equation}
does not vanish on the set $A$ introduced in (\ref{aa4}).
Since $p$ is a polynomial with real nonpositive zeros only,
the polynomial $q$ can be written as follows
\be
\label{aa17}
q(z)=q_{0}\prod_{j=1}^{m}(q_{j}+z^{2}),\ \ m = {\rm deg}p , \ \  q_{j}\geq 0.
\ee
It is no need to consider the trivial case of constant $p$.
Consider the simplest nontrivial case where all $q_j = 0, \ \ j=1, \ \dots , m$
and $q_0 \neq 0$,
that is
$$
q(z) = q_0 z^{2m} , \ \ m\geq 1 .
$$
Then one gets
$$
(\k + \Lambda_\th )q (z) = [\k z^2 + m( \th +m - 1 )]q_0 z^{2(m-1)} ,
$$
that obviously does not vanish on $A$. From now on we suppose
that in the product in (\ref{aa17}) there is
at least one positive $q_j$. Then one may write
$$
Dq(z)=2zq(z)r(z),
$$
where
$$
r(z)=\sum_{j=1}^{m}\frac{1}{q_{j}+z^{2} }.
$$
Then
\beq
(\k +\Lambda _{\theta })q(z) & = & \left[ \k +
\left( \theta +\frac{z}{2}%
D\right)\left( \frac{1}{2z}D\right)
\right] q(z) \nonumber \\
& = & \k q(z)+\theta q(z) r(z)+\frac{z}{2} D [q(z) r(z)]. \nonumber
\eeq
We set
\begin{equation}
\label{a18}
Q(z) \ \sta \  q(z) (\k +\theta r(z) ),\quad
P(z) \ \sta \ \frac{z}{2},\quad
Q_{1}(z) \ \sta \ q(z) r(z).
\end{equation}
Let us show that
\begin{equation}
\label{a19}
R(v,w) \ \sta \ Q(w)+P(w)vQ_{1}(w)\neq 0,
\end{equation}
whenever $v,$ $w\in A.$ To this end we rewrite the latter
$$
R(v,w)=\half (2\theta +v w)q(w)\left[ \frac{2\k }{2\theta +v w}+%
r(w)\right] .
$$
Obviously $(2\theta +v w)\neq 0$ whenever $v,$ $w\in A.$
The same property possesses also $q(w) = p(w^2 )$, $p\in {\cal P}^+ $.
Therefore, the eventual vanishing of $R(v , w)$ would imply
$$
\frac{2\k }{2\theta +v w}+\sum_{j=1}^{m}\frac{1}{q_{j}+w^{2}}=0
$$
or equivalently
\begin{equation}
\label{a20}
\frac{2\k (2\theta +\bar{v}\bar{w})}{\left| 2\theta +vw\right|^{2}}
+\sum_{j=1}^{m}\frac{q_{j}+\bar{w}^{2}}{\left|
q_{j}+w^{2}\right|^{2}}=0.
\end{equation}
Introduce
$$
A(v,w) \ \sta \ \frac{4\k \theta }{\left| 2\theta +v w\right|^{2}}%
+\sum_{j=1}^{m}\frac{q_{j}}{\left| q_{j}+w^{2}\right|^{2}}>0
$$
$$
B(v,w) \ \sta \ \frac{2\k }{\left| 2\theta +v w \right|^{2} }\geq 0,\ \
\ C(v,w)\ \sta \ \sum_{j=1}^{m}\frac{1}{\left| q_{j}
+w^{2}\right|^{2} }>0
$$
and rewrite (\ref{a20}) as follows
$$
C\bar{w}^{2}+B\bar{v}\bar{w}+A=0.
$$
Thus
$$
\bar{v}=-\frac{A}{B \left| w \right|^{2}} w-\frac{C}{B} \bar{w},\ \ \ B>0 ,
$$
or
$$
\bar{w}^{2}=-\frac{A}{C}<0,\ \ \ B=0.
$$
The above expressions imply the following conclusions
$$
{\rm Re}v=-\left( \frac{A}{B \left| w\right|^{2} }+\frac{C}{B} \right)
 {\rm Re}w<0 ,
$$
or
$$
{\rm Re}w=0.
$$
Both ones run in counter with the assumption that $v$ and $w$ belong to $A$.
This
means $R(v,w)\neq 0,$ whenever $v,$ $w\in A.$ By means of Lemma
\ref{a3lm} one concludes that $(\k +\Lambda _{\theta })q(z)\neq 0,$ for
$z\in A$.
\kasten

 {\bf Proof of Theorem \ref{a1tm}.} It suffices to show that for
any $%
\varphi \in {\cal L}_{a}^{+}$ and $f\in {\cal L}_{b}^{+},$  the function
 $\varphi (\Delta _{\theta })f$ belongs to $ {\cal L}^{+} .$
To this end we choose a sequence
$\{\varphi _{n} ,  \ n\in \N \}\subset {\cal P}^{+} $, converging to $\varphi $ in $%
{\AC} _{a}$ (see Corollary \ref{n11co}).
The fact $\vp_n \in {\cal P}^{+}$
implies
$$
\varphi _{n}(z)=\phi_{n}\prod_{j=1}^{m_{n}}(\k _{j,n}+z),\ \ \
\phi_{n} \neq 0,\quad \k _{j,n}\geq 0.
$$
Then
$$
g_{n}(z)\ \sta \ \varphi _{n}(\Delta _{\theta
})f(z)=\phi_{n}\prod_{j=1}^{m_{n}}(\k _{j,n}+\Delta _{\theta })f(z).
$$
Thus Lemma \ref{a4lm} yields $g_{n}(z)\in {\cal L}^{+}.$ Corollary \ref {1co}
implies that $\{g_{n}\}$ converges to $\varphi (\Delta _{\theta })f$ in
${\AC}_{c}.$ The latter yields $\varphi (\Delta
_{\theta })f\in \cal{L}^{+}$.
\kasten

{\bf Proof of Theorem \ref{e1tm}.}
Consider the operator valued function $ (0,t_0 )\ni
t \mapsto \vp_t (\D_\th )
\in {\bf B}({\cal B}_b , {\cal B}_c )$, where the latter is
the Banach space of all linear bounded operators between the Banach
spaces ${\cal B}_b$, ${\cal B}_c$ with $c =b(1-t_0 b)$, $t_0 b<1$.
One may easily show that this function is continuous and differentiable
in the norm-topology and its derivative is
$$
\vp'_t(\Delta_\theta ) = \Delta_\theta \vp_t (\Delta_\theta ).
$$
Therefore, for $t\in (0, 1/\ve )$, one has
$$
\frac{\pp f(t, z)}{\pp t} = \Delta_\theta \left( \exp (t\Delta_\theta  )
g \right) (z),
$$
which proves the first line in (\ref{e2}). The second line,
which gives the extension of the solution to all positive values
of $t$, is easily obtained from the first one by means of
Proposition \ref{a2pn}. Further, we substitute in (\ref{e2}) the
initial condition in the form (\ref{e3}) and apply Proposition
\ref{a3pn}. This yields
\beq
\label{b1}
f(t,z) & = & (1+\ve t)^{-\theta}\exp \left(-\frac{\ve z}{1+ \ve t}
\right)
\left[ \exp\left(\frac{t}{1+ \ve t}\Delta_\theta  \right) h \right]
\left(\frac{z}{(1+\ve t)^2} \right) \nonumber\\
& \sta &  (1+\ve t)^{-\theta}\exp \left(-\frac{\ve z}{1+ \ve t}\right)
h_t \left(\frac{z}{(1+\ve t)^2}\right).
\eeq
By Corollary \ref{1co}
$h_t \in {\AC}_0 $, and by Theorem \ref{a2tm} $h_t \in {\LC}_0 $
if $h \in {\LC}_0 $.
The former proves claim (i) and the latter does claim (iii).
It remains to prove the convergence stated in (ii). The mentioned
continuity of the operator $\exp (t \Delta_\theta )$ implies that
in ${\AC}_0$
\be
\label{b2}
\left[ \exp\left(\frac{t}{1+ \ve t}\Delta_\theta  \right) h \right]
\left(\frac{z}{(1+\ve t)^2} \right)
\ \ra \
\left\{ \exp \left(\frac{1}{\ve }\Delta_\theta  \right)h \right\}(0).
\ee
Therefore, the product in (\ref{b1}) tends to zero in ${\AC}_\ve$
when $t\ra +\infty $.
\kasten
{\bf Proof of the identity (\ref{aa23}).}
Assume that $k\leq m $, then the left hand side of (\ref{aa23})
may be brought into the following form
$$
\frac{1}{\Gamma (z+k)}\prod_{l=1}^{k}(z+m +k-l),
$$
which one rewrites as
$$
\frac{1}{\Gamma (z+k)}\prod_{l=1}^{k}[(z+k-j_l) +(m+j_l -l)],
$$
with arbitrary $j_l$. Then one opens the brackets $[ \ . \ ]$ and
transforms the product of sums into the sum of the products
choosing in every term an appropriate $j_l \in \{ 1, \dots , l\}$.
This yields
\beq
& & \frac{1}{\Gamma (z+k)}  \sum_{n=0}^{k}
\left[\prod_{l=1}^{k-n}[(z+k-l)m(m-1)\dots (m-n+1){{k}\choose{n}}
\right] \nonumber \\
&= & \sum_{n=0}^{k}n!{{m}\choose{n}}{{k}\choose{n}}
\frac{1}{\Gamma (z+k)}
\prod_{l=1}^{k-n}(z+k-l) \nonumber \\
& =& \sum_{n=0}^{k}{{m}\choose{n}}{{k}\choose{n}}
\frac{n!}{\Gamma (z+k)} = {\rm RHS}(\ref{aa23}).
\eeq
\kasten

{\bf Acknowledgment} This work was brought into the final form
during the stay of Yuri Kozitsky in Bochum financially supported
by SFB--237 (Essen-Bochum-D{\"u}sseldorf), which is gratefully
acknowledged.

\vskip.6cm

      Yuri Kozitsky \\[.5cm]
        Institute of Mathematics,\\
                Maria Curie-Sklodowska University\\
               PL 20-031 Lublin (Poland)\\ [.2cm]
      Institute for Condensed Matter Physics,\\
      Lviv (Ukraine)\\[.2cm]
       {\sf e-mail jkozi@golem.umcs.lublin.pl}\\
\vskip.6cm

Lech Wolowski \\[.5cm]
Institute of Mathematics,\\
                Maria Curie-Sklodowska University\\
               PL 20-031 Lublin (Poland);\\ [.2cm]
       {\sf e-mail lechw@golem.umcs.lublin.pl}

\end{document}